\documentclass[12pt]{amsart}
\usepackage{amsmath,amssymb,enumerate}  % mathtools
\usepackage[margin=1in]{geometry}

\usepackage{amsmath}
\usepackage[subnum]{cases}

\usepackage{tikz}
\usetikzlibrary{arrows}
\usetikzlibrary{positioning}

\newtheorem{theorem}{Theorem}
\newtheorem{lemma}{Lemma}

\DeclareMathOperator{\mex}{mex}

\def\gl{\lambda}

\def\U{\mathcal{U}}
\def\G{\mathcal{G}}
\def\T{\mathcal{T}}

\def\ZZP{\mathbb{Z}_{+}}
\def\ZZ{\mathbb{Z}}

\def\cH{\mathcal{H}}

\begin{document}

\title{Sprague-Grundy Function of Symmetric Hypergraphs}

%%%%%%%%%%%%%%%%%%%%%%%%%%%%%%%%%%%%%%%%%%%%%%%%%%%%%%%
\author{Endre Boros}
\address{MSIS and RUTCOR, RBS, Rutgers University,
100 Rockafeller Road, Piscataway, NJ 08854}
\email{endre.boros@rutgers.edu}

\author{Vladimir Gurvich}
\address{MSIS and RUTCOR, RBS, Rutgers University,
100 Rockafeller Road, Piscataway, NJ 08854; \\
Dep. of Computer Sciences, National Research University,
Higher School of Economics (HSE), Moscow}
\email{vladimir.gurvich@rutgers.edu}

\author{Nhan Bao Ho}
\address{Department of Mathematics and Statistics, La Trobe University, Melbourne, Australia 3086}
\email{nhan.ho@latrobe.edu.au, nhanbaoho@gmail.com}

\author{Kazuhisa Makino}
\address{Research Institute for Mathematical Sciences (RIMS)
Kyoto University, Kyoto 606-8502, Japan}
\email{e-mail:makino@kurims.kyoto-u.ac.jp}

\author{Peter Mursic}
\address{MSIS and RUTCOR, RBS, Rutgers University,
100 Rockafeller Road, Piscataway, NJ 08854}
\email{peter.mursic@rutgers.edu}

%\authornames{Endre Boros \and Vladimir Gurvich \and Nhan Bao Ho \and Kazuhisa Makino \and Peter Mursic}

%%%%%%%%%%%%%%%%%%%%%%%%%%%%%%%%%%%%%%%%%%%%%%%%%%%%%%%

\subjclass[2000]{91A46}

\keywords{Impartial game, Sprague-Grundy function, $NIM$, hypergraph $NIM$, JM hypergraph, symmetric hypergraph, transversal-free hypergraph.}

\begin{abstract}
We consider a generalization of the classical
game of $NIM$ called hypergraph $NIM$.
Given a hypergraph  $\cH$
on the ground set  $V = \{1, \ldots, n\}$ of $n$ piles of stones,
two players alternate in choosing a hyperedge
$H \in \cH$ and strictly decreasing all piles $i\in H$.
The player who makes the last move is the winner.
Recently it was shown that for many classes of hypergraphs the Sprague-Grundy function of the corresponding game is given by the formula introduced originally by Jenkyns and Mayberry (1980).
In this paper we characterize symmetric hypergraphs for which the Sprague-Grundy function is described by the same formula.
\end{abstract}

\maketitle

\section{Introduction}
\label{s1}

In the classical game of $NIM$ there are
$n$  piles of stones and two players move alternating.
A move consists of choosing a nonempty pile and taking some positive number of stones from it.
The player who cannot move is the looser.
Bouton \cite{Bou901} analyzed this game and described the winning strategy for it.

In this paper we consider the following generalization of $NIM$.
Given a hypergraph $\cH \subseteq 2^V$, where $V = \{1, \dots , n\}$, two players alternate in choosing
a hyperedge $H \in \cH$ and strictly decreasing all piles $i\in H$.
We assume in this paper that $\cH\not=\emptyset$, $\emptyset \not\in \cH$, and $|V|=n$ for all considered hypergraphs $\cH\subseteq 2^V$.
In other words, every move strictly decreases some of the piles.
Similarly to $NIM$, the player who cannot move is loosing.
This game is denoted by $NIM_\cH$. Such games were introduced in \cite{BGHMM15,BGHMM16,BGHMM18} and called hypergraph $NIM$.

Hypergraph $NIM$ games are impartial.
In this paper we do not need to immerse in the theory of impartial games.
We will need to recall only a few basic facts to explain and motivate our research.
We refer the reader to \cite{Alb07,BCG01-04} for more details.

It is known that the set of positions of an impartial game can
uniquely be partitioned into sets of winning and loosing positions.
Every move from a loosing position goes to a winning one, while
from a winning position there always exists a move to a loosing one.
This partition shows how to win the game, whenever possible.
The so-called Sprague-Grundy (SG) function $\G_\Gamma$ of an impartial game $\Gamma$ is
a refinement of the above partition, see Section \ref{s2} for the precise definition.
Namely, $\G_\Gamma(x)=0$ if and only if $x$ is a loosing position.
The notion of the SG function for impartial games was introduced by Sprague and Grundy
\cite{Spr35, Spr37, Gru39} and it plays a fundamental role in determining the winning-loosing partition
of \emph{disjunctive sums} of impartial games.

Finding a formula for the SG function of an impartial game remains a challenge.
Closed form descriptions are known only for some special classes of impartial games.
We recall below some known results.
The purpose of our research is to extend these results and
to describe classes of hypergraphs for which
we can provide a closed formula for the SG function of $NIM_\cH$.

The game $NIM_\cH$ is a common generalization of several families of impartial games considered in the literature. For a subset $S\subseteq V$ and integer $1\leq k\leq |S|$ we denote by
\[
\binom{S}{k}=\{H\subseteq S\mid |H|=k\}.
\]
For instance, if $\cH=\binom{V}{1}$ then $NIM_\cH$ is the classical $NIM$.
The case of $\cH=\bigcup_{j=1}^k \binom{V}{k}$, where $k<n$, was considered by Moore \cite{Moo910}.
He characterized for these games the set of loosing positions, that is those with SG value $0$.
Jenkyns and Mayberry \cite{JM80} described also the
set of positions in which the SG value is $1$ and
provided an explicit formula for the SG function in case of $k = n-1$.
This result was extended in \cite{BGHM15}.
In \cite{BGHMM15} the game $NIM_\cH$ was considered
in the case of $\cH=\binom{V}{k}$
and the corresponding SG function was determined when  $2k\geq n$. Further examples for hypergraph $NIM$ such as matroid and graph $NIM$  were given in \cite{BGHMM18}. Surprisingly, for many of these examples the SG function is described by the same formula, special case of which was introduced by Jenkyns and Mayberry \cite{JM80}. In honor of their contribution, this formula and the hypergraphs for which it describes the SG function were called \emph{JM} in \cite{BGHMM18}.

In this paper we characterize symmetric JM hypergraphs. A hypergraph $\cH\subseteq 2^V$ is called \emph{symmetric} if $H\in \cH$ implies that all subsets of $V$ of size $|H|$ are also hyperedges of $\cH$. Note that many of the above cited examples are symmetric hypergraphs.

To state our main result formally we need to introduce some additional notation.
We denote by $\ZZP$ the set of nonnegative integers and
use $x\in\ZZP^V$ to describe a position, where coordinate $x_i$ denotes
the number of stones in pile $i \in V$.
Given a hypergraph $\cH$ and position $x \in \ZZP^V$, we denote by $\G_\cH(x)$ the SG value of $x$ in $NIM_\cH$.
The Tetris value $\T_{\cH}(x)$ was defined in \cite{BGHMM15} as the maximum number of consecutive moves
that the players can make in $NIM_\cH$ starting from position $x$.

To a position $x \in \ZZP^V$ of $NIM_\cH$ let us associate the following quantities:
\begin{subequations}\label{e-myv}
\begin{align}
m(x) &=\min_{i\in V} x_i \label{e-m}\\
y^{}_\cH(x) &=\T_\cH(x-m(x)e)+1 \label{e-y}\\
v^{}_\cH(x) &=\binom{y^{}_\cH(x)}{2} +
\left( \left(m(x)-\binom{y^{}_\cH(x)}{2}-1\right)\mod y^{}_\cH(x)\right),\label{e-v}
\end{align}
\end{subequations}
where $e$ is the $n$-vector of full ones. Finally, we define
%\begin{subequations}\label{e-JM}
%\begin{align}[left ={\U_\cH(x) ~=~ \empheqlbrace}]
%\T_\cH(x) & ~~\text{ if } m(x) \leq \binom{y^{}_\cH(x)}{2}\label{e-JM-I}\\
%v^{}_\cH(x) & ~~\text{ if } m(x) > \binom{y^{}_\cH(x)}{2}.\label{e-JM-II}
%\end{align}
%\end{subequations}
\begin{numcases}{\U_\cH(x) =}
\T_\cH(x)   \quad \text{ if } m(x) \leq \binom{y^{}_\cH(x)}{2}\label{e-JM-I}\\
v^{}_\cH(x) \quad \text{ otherwise.} \label{e-JM-II} %m(x) > \binom{y^{}_\cH(x)}{2}.\label{e-JM-II}
\end{numcases}

The expression \eqref{e-JM-I}-\eqref{e-JM-II} is
the \emph{JM formula} mentioned earlier.

Given a hypergraph $\cH\subseteq 2^V$ and a subset $S\subseteq V$, we denote by $\cH_S$ the \emph{induced subhypergraph}, defined as
\[
\cH_S=\{H\in\cH\mid H\subseteq S\}.
\]
A set $T \subseteq V$ is called a {\em transversal} if  $T \cap H\neq\emptyset$ for all $H\in\cH$.
A hypergraph $\cH$ is called {\em transversal-free} if no hyperedge $H \in \cH$ is a transvesal of $\cH$.
Finally, we say that $\cH$ is \emph{minimal transversal-free} if it is transversal-free and
every nonempty proper induced subhypergraph of it is not.

Consider an integer sequence $\gl=(\gl_1,\gl_2,\ldots,\gl_k)$ such that $0<\gl_1<\gl_2<\cdots <\gl_k\leq n$ and an associated hypergraph
\begin{equation}\label{e-spectrum}
\cH(\gl) ~=~ \bigcup_{j=1}^k \binom{V}{\gl_j}.
\end{equation}
The sequence $\gl$ is called the \emph{spectrum} of the symmetric hypergraph $\cH(\gl)$.
All symmetric hypergraphs have a spectrum and arise in this way (recall that a hypergraph  considered in this paper is not $\emptyset$ and does not contain $\emptyset$).

The following two theorems summarize our main results.

\begin{theorem}\label{t-1}
A symmetric hypergraph is JM if and only if it is minimal transversal-free and $n\geq 3$.
\end{theorem}

\begin{theorem}\label{t-2}
A symmetric hypergraph $\cH(\gl)$ defined by \eqref{e-spectrum} is minimal transversal-free if and only if its spectrum satisfies the following relations:
\begin{itemize}
\item[(i)] $\gl_{i+1} -\gl_{i} \leq \gl_1$ for all $i=1,\ldots,k-1$.
\item[(ii)] $\gl_1+\gl_k=n$.
\end{itemize}
\end{theorem}

\smallskip

Theorems \ref{t-1} and \ref{t-2} extend two previous results stating that
$\bigcup_{j=1}^{n-1}\binom{V}{k}$ and $\binom{V}{n/2}$ are JM hypergraphs \cite{BGHMM15,JM80}.

Let us also add that by Theorem 1 in \cite{BGHMM18} we can derive numerous non-symmetric JM hypergraphs from the above family of symmetric ones.

We finally remark that it is not easy to describe the closed form of the SG  function of symmetric hypergraph $NIM$ games, in general.
For example, the case
$\binom{V}{1} \cup \binom{V}{2}$ for $n=4$ seems to be difficult. In this case at least the loosing positions are know, due to \cite{Moo910}. For the case of
$\binom{V}{2}$ for $n=5$ we are not even aware of a useful characterization of the set of loosing positions. In contrast, $\binom{V}{2}\cup\binom{V}{3}$ is a symmetric JM hypergraph by our theorems above.

\smallskip

The structure of the paper is as follows. In Section \ref{s2}, we fix our notation and define basic concepts.
We also provide several properties of JM hypergraphs that are needed for our proofs later.
In Subsections \ref{s-pt1} and \ref{s-pt2} we prove Theorems \ref{t-1} and \ref{t-2}, respectively.

%Finally, in Section \ref{sec-cr}, we present further examples
%of JM hypergraphs and discuss related topics.

\section{Basic Concepts and Notation}
\label{s2}

In this section we introduce the basic notation and definitions. We prove some of the basic properties of $NIM_\cH$ games, the Tetris functions and the JM formula.

\subsection{Hypergraph $NIM$ Games}
\label{ss21}

We need to recall first the precise definition of impartial games and the SG function.

To a subset $S\subseteq \ZZP$ of nonnegative integers let us associate its \emph{minimal excludant}
$\mex(S)=\min\{i\in \ZZP\mid i\not\in S\}$, that is,
the smallest nonnegative integer that is not included in $S$. Note that $\mex(\emptyset)=0$.

An impartial game $\Gamma$
is played by two players over a (possibly infinite) set $X$ of positions. They take turns
to move, and the one who cannot move is the looser. For a position $x\in X$ let us denote by $N^+(x)\subseteq X$
the set of positions $y\in N^+(x)$ that are reachable from $x$ by a single move.
For $y\in N^+(x)$ we denote by $x\to y$ such a move. We assume that the same set of moves are available for both players from every position.
We also assume that no matter how the players play and which position they start, the game ends in a finite number of moves.
The SG function $\G_\Gamma$ of the game is a mapping $\G_\Gamma:X\mapsto \ZZP$ that associates a nonnegative integer to every position, defined by the following  recursive formula:
\[
\G_\Gamma(x) ~=~ \mex\{\G(y)\mid y\in X, ~s.t.~ \exists x\to y\} .
\]

In our proofs we shall use the following, more combinatorial characterization of SG functions that can be derived easily from the above definition.
Assume that $\Gamma$ is an impartial game over the set of positions $X$, and
$g:X\to \ZZP$ is a given function. Then, $g$
is the SG function of $\Gamma$ if and only if the following two conditions hold:
\begin{itemize}
\item For all positions $x\in X$ and moves $x\to y$ in $\Gamma$ we have $g(x)\neq g(y)$.
\item For all positions $x\in X$ and integers $0\leq z<g(x)$ there exists a move $x\to y$ in $\Gamma$ such that $g(y)=z$.
\end{itemize}

It is easy now to verify that for a hypergraph $\cH\subseteq 2^V$ the game $NIM_\cH$ is indeed an impartial game over the infinite set $X=\ZZP^V$ of positions.

\medskip

Let us note that all quantities used in \eqref{e-m}--\eqref{e-JM-II}, namely $m(x)$, $y_\cH(x)$, $v_\cH(x)$, $\T_\cH(x)$ as well as $\U_\cH$ are well defined
for an arbitrary hypergraph $\cH$.  Let us also note that the values $m(x)$ and $y_\cH(x)$ determine completely the value of $v_\cH(x)$.

To simplify our language and notation in the sequel, let us call a position $x\in \ZZP^V$
\emph{long} (in $NIM_\cH$) if $m(x) \leq \binom{y^{}_\cH(x)}{2}$
(that is, if  $\U_\cH(x)=\T_\cH(x)$) and call it
\emph{short} if $m(x) > \binom{y^{}_\cH(x)}{2}$ (that is, if $\U_\cH(x)=v_\cH(x)$).

Let us first recall from \cite{BGHMM18} some essential necessary and sufficient conditions for a hypergraph to be JM.

\begin{lemma}[\cite{BGHMM18}]\label{l3}
If $\cH$ is a JM hypergraph, then it is minimal transversal-free.
\end{lemma}

\begin{lemma}[\cite{BGHMM18}]\label{l-ns-JM}
A hypergraph $\cH\subseteq 2^V$ is JM if the following four conditions hold:
\begin{itemize}
\item[(A0)] $\cH$ is transversal-free.
\item[(B1)] For every long %%type I
position $x\in\ZZP^V$ and integer $m(x)\leq z <\T_\cH(x)$ there exists a move $x\to x'$ such that $x'$ is long %%of type I
and $\T_\cH(x')=z$.
\item[(C2)] For every position $x\in\ZZP^V$ and integer $1\leq \eta < y_\cH(x)$ there exists a move $x\to x'$ such that $m(x')=m(x)$ and $y_\cH(x')=\eta$.
\item[(C3)] For every position $x\in\ZZP^V$ and integers $0\leq \mu < m(x)$ and $m(x)-\mu +1 \leq \eta \leq y_\cH(x)$ there exists a move $x\to x'$ such that $m(x')=\mu$ and $y_\cH(x')=\eta$.
\end{itemize}
\end{lemma}

We shall also need to state some properties of the Tetris function, in particular of the Tetris function of symmetric hypergraphs.

According to the rules of $NIM_\cH$, if $x\to x'$ is a move then for the set $H=\{i\mid i\in V,~ x_i>x'_i\}$ we must have $H\in \cH$. We call such a move an $H$-\emph{move}.
For a subset $S\subseteq V$ let us denote by $\chi(S)$ the characteristic vector of $S$, that is, $\chi(S)_j=1$ if $j\in S$ and $\chi(S)_j=0$ if $j\not\in S$.
We denote for a position $x\in \ZZP^V$ and hyperedge $H\in\cH$ by $x^{s(H)}$ the vector $x-\chi(H)$. Thus, $x\to x^{s(H)}$ is an $H$-move in $NIM_\cH$. We call such a move also a \emph{slow move}, since we have $x'\leq x^{s(H)}$ for all $H$-moves $x\to x'$. Note that in the definition of the Tetris function, it is enough to consider slow moves using only inclusion-wise minimal hyperedges.

Let us first recall a few basic properties of Tetris functions of hypergraph $NIM$ games. The proof is simple and follows by the definition of the Tetris function, see  \cite{BGHMM18}.

\begin{lemma}\label{l-t-contiguity}
For an arbitrary hypergraph $\cH\subseteq 2^V$ and positions $x',x''\in\ZZ_+^V$ the inequality $x''\leq x'$ implies $\T_\cH(x'')\leq \T_\cH(x')$ and we have
\[
\{\T_\cH(y) \mid x''\leq y\leq x'\} ~=~ \{ t \mid \T_\cH(x'')\leq t\leq \T_\cH(x')\}.
\]
Furthermore, if for a position $x\in\ZZ_+^V$ we have both $x\to x'$ and $x\to x''$ as $H$-moves for some $H\in \cH$, then all moves $x\to y$ for $x''\leq y\leq x'$ are $H$-moves.
\end{lemma}

Let us next show some additional properties that of Tetris functions of more special hypergraphs.

Let us next recall a lemma from \cite{BGHM15}
\begin{lemma}\label{l-s2}
Let $\cH(\gl)$ be a symmetric hypergraph as defined in \eqref{e-spectrum}. Then for any position $x\in\ZZP^V$ and indices $i\neq j$ such that $x_i>x_j$ we have $\T_k(x')\geq \T_k(x)$, where $x'\in\ZZP^V$ is defined by
\[
x'_\ell = \begin{cases} x_i-1 & \text{ if } \ell =i,\\ x_j+1 & \text{ if } \ell=j,\\ x_\ell & \text{ otherwise}.\end{cases}
\]
\qed
\end{lemma}

\begin{lemma}\label{l-l1}
Given a symmetric hypergraph $\cH=\cH(\gl)\subseteq 2^V$ and a position $x\in\ZZ_+^V$ such that $x_1\geq x_2\geq \cdots \geq x_n>0$, let us consider the following two hyperedges of size $\gl_1$ each: $H_1=\{1,2,\dots ,\gl_1\}$ and $H_2=(H_1\setminus\{1\})\cup\{ n\}$. Then we have $\T_\cH(x^{s(H_i)}) = \T_\cH(x)-1$ for both $i=1,2$.
\end{lemma}

\proof
Let us first consider an arbitrary hyperedge $H\in\cH$ of size $\gl_1$ such that $\T_\cH(x^{s(H)})=\T_\cH(x)-1$. By the definition of the Tetris function such a hyperedge exists since $\T_\cH(x)>0$. Then we can apply Lemma \ref{l-s2} repeatedly and conclude that
$\T_\cH(x^{s(H_1)})= \T_\cH(x^{s(H)})=\T_\cH(x)-1$.

Let us now consider a longest sequence of consecutive slow moves with hyperedges $H^1$, $H^2$, \dots , $H^{\T_\cH(x)}$ where $H^1=H_1$. Assume first that there exists a hyperedge $H^i\ni n$. Then $\T_\cH(x^{s(H^i)})=\T_\cH(x)-1$ and thus we can apply Lemma \ref{l-s2} and conclude that $\T_\cH(x^{s(H_2)})= \T_\cH(x^{s(H^i)})=\T_\cH(x)-1$. Otherwise, if no such hyperedge exists, then $H_2$, $H^2$, \dots , $H^{\T_\cH(x)}$ also forms a longest sequence of consecutive slow moves, proving our claim.
\qed

\section{Proof of the Main Results}\label{s-p}

\subsection{Proof of Theorem \ref{t-2}}\label{s-pt2}

Let us assume first that $\cH=\cH(\gl)$ is minimal transversal-free.

If there is an index $i$ such that $\gl_{i+1}-\gl_i>\gl_1$, then let us consider a proper subset $S\subseteq V$ of size $|S|=\gl_{i+1}-1<n$. The induced subhypergraph $\cH_S$ has no transversal hyperedge, contradicting our assumption. Thus we must have $\gl_{i+1}-\gl_i\leq \gl_1$ for all indices $i=1,\dots ,k-1$.

If $\gl_1+\gl_k < n$, then consider a subset $S\subseteq V$ of size $|S|=n-1$. It is easy to see again that $\cH_S$ has no transversal hyperedge, leading to a contradiction, as above.
Finally, if $\gl_1+\gl_k > n$, then $\cH$ has a transversal hyperedge, contradicting our assumption. Thus, we must have $\gl_1+\gl_k = n$.

Assume next that conditions (i) and (ii) of Theorem \ref{t-2} hold. Then condition (ii) implies that $\cH$ has no transversal hyperedge. Let us consider an arbitrary proper subset $\emptyset\neq S\subsetneq V$. If $|S| < \gl_1$, then $\cH_S$ is an empty hypergraph. If $\gl_i \leq |S| < \gl_{i+1}$ for some index $1\leq i \leq k$ (assuming $\gl_{k+1}=n$), then any hyperedge of size $\gl_i$ inside $S$ is a transversal of $\cH_S$. Therefore conditions (i) and (ii) imply that $\cH=\cH(\gl)$ is minimal transversal-free.
\qed

\subsection{Proof of Theorem \ref{t-1}}\label{s-pt1}

It is easy to verify that if $n\leq 2$, then there exist no JM hypergraphs. Thus, we can assume in the sequel that $n\geq 3$.

Observe next that by Lemma \ref{l3} a JM hypergraph must be minimal transversal-free.

For the reverse direction we consider a symmetric minimal transversal-free hypergraph $\cH(\gl)$. If $\gl_1=1$ then by Theorem \ref{t-2} we have $\gl=(1,2,\dots ,n-1)$. Consequently the hypergraph $\cH(\gl)$ coincides with the one considered in \cite{JM80}, and thus their result implies our claim.

For the remaining cases, when $\gl_1>1$, we show that the sufficient conditions of Lemma \ref{l-ns-JM} hold. We break this proof into three technical lemmas, and start with the simplest one. For positions $a,b\in\ZZ_+^V$, $a\leq b$ we define $[a,b]=\{x\in\ZZ_+^V\mid a\leq x\leq b\}$.

\begin{lemma}\label{l-C2}
If $\cH=\cH(\gl)$ is a symmetric hypergraph such that its spectrum $\gl$ satisfies conditions (i) and (ii) of Theorem \ref{t-2}, then condition (C2) holds.
\end{lemma}

\proof
Let us consider a position $x\in\ZZ_+^V$ such that $x_1\geq x_2\geq \cdots \geq x_n$ and $y_\cH(x)>1$.
Let us define $\ell=\max\{j\mid x_{\gl_j}>m(x)\}$. By the assumption $y_\cH(x)>1$ it is well defined.

Next we define $\gl_0=0$ and positions $a^i,b^i$ for $i=1,\dots ,\ell$ as follows.

\[
a^i_j = \begin{cases}
m(x) & \text{if } j\leq \gl_{i-1},\\
x_j-1 & \text{if } \gl_{i-1}<j\leq \gl_{i},\\
x_j & \text{otherwise},
\end{cases}
~~~~~~~
b^i_j = \begin{cases}
m(x) & \text{if } j\leq \gl_{i},\\
x_j & \text{otherwise}.
\end{cases}
\]

Note first that for all indices $i=1,\dots ,\ell$ we have $a^i\geq b^i$, and there exists a hyperedge $H^i\in\cH$ such that both $x\to a^i$ and $x\to b^i$ are $H^i$-moves.

Let us note next that due to the definition of $y_\cH$, condition (i) of Theorem \ref{t-2}, and Lemma \ref{l-l1} we have $y_\cH(a^{i+1}) \geq y_\cH(b^i)-1$ for $i=1,\dots ,\ell-1$. Furthermore $y_\cH(a^1)= y_\cH(x)-1$ by Lemma \ref{l-l1} and $y_\cH(b^\ell)=1$.

Let us then note that we have $m(x')=m(x)$ for all $x'\in \bigcup_{i=1}^\ell [a^i,b^i]$. Furthermore we can apply Lemma \ref{l-t-contiguity} to the pairs $(a^i,b^i)$, $i=1,\dots ,\ell$ and obtain
\[
\bigcup_{i=1}^\ell \{y_\cH(x')\mid x'\in [a^i,b^i]\} ~=~ [1,y_\cH(x)-1].
\]
\qed

\begin{lemma}\label{l-C3}
If $\cH=\cH(\gl)$ is a symmetric hypergraph such that its spectrum $\gl$ satisfies conditions (i) and (ii) of Theorem \ref{t-2} and $\gl_1>1$, then condition (C3) holds.
\end{lemma}

\proof
Let us consider a position $x\in\ZZ_+^V$ such that $x_1\geq x_2\geq \cdots \geq x_n$ and $y_\cH(x)>1$.
Next we define positions $a^i,b^i$ for $i=0,\dots ,k$ as follows.

\[
a^0_j = \begin{cases}
x_j-1 & \text{if } j< \gl_{1},\\
\mu & \text{if } j=n,\\
x_j & \text{otherwise},
\end{cases}
~~~~~~~
b^0_j = \begin{cases}
\mu & \text{if } j< \gl_{1} \text{ or } j=n,\\
x_j & \text{otherwise}.
\end{cases}
\]

\[
a^1_j = \begin{cases}
\mu & \text{if } j< \gl_{1},\\
x_j-1 & \text{if } j=\gl_1,\\
x_j & \text{otherwise},
\end{cases}
~~~~~~~
b^1_j = \begin{cases}
\mu & \text{if } j\leq  \gl_{1},\\
x_j & \text{otherwise}.
\end{cases}
\]
Note that $m(a^1)=\mu$ because we assumed $\gl_1>1$.

For $i=2,\dots ,k$ we set

\[
a^i_j = \begin{cases}
\mu & \text{if } j\leq \gl_{i-1},\\
x_j-1 & \text{if } \gl_{i-1}<j\leq \gl_{i},\\
x_j & \text{otherwise},
\end{cases}
~~~~~~~
b^i_j = \begin{cases}
\mu & \text{if } j\leq \gl_{i},\\
x_j & \text{otherwise}.
\end{cases}
\]

Note first that for all indices $i=0,\dots ,k$ we have $a^i\geq b^i$, and there exists a hyperedge $H^i\in\cH$ such that both $x\to a^i$ and $x\to b^i$ are $H^i$-moves.

Let us note next that due to the definition of $y_\cH$, condition (i) of Theorem \ref{t-2}, and Lemma \ref{l-l1} we have $y_\cH(a^{i+1}) \geq y_\cH(b^i)-1$ for $i=1,\dots ,k-1$. Furthermore $y_\cH(a^0)\geq y_\cH(x)$ and $y_\cH(b^k)=m(x)-\mu+1$ because of condition (ii) of Theorem \ref{t-2}.

Let us then note that we have $m(x')=\mu$ for all $x'\in \bigcup_{i=0}^k [a^i,b^i]$ because we assumed $\gl_1>1$. Furthermore we can apply Lemma \ref{l-t-contiguity} to the pairs $(a^i,b^i)$, $i=0,\dots ,k$ and obtain
\[
\bigcup_{i=0}^k \{y_\cH(x')\mid x'\in [a^i,b^i]\} \supseteq [m(x)-\mu+1,y_\cH(x)].
\]
\qed

\begin{lemma}\label{l-B1}
If $\cH=\cH(\gl)$ is a symmetric hypergraph such that its spectrum $\gl$ satisfies conditions (i) and (ii) of Theorem \ref{t-2} and $\gl_1>1$, then condition (B1) holds.
\end{lemma}

\proof
Let us consider a long position $x\in\ZZ_+^V$ such that $x_1\geq x_2\geq \cdots \geq x_n$.  Next we define positions $a^i,b^i$ for $i=0,\dots ,k$ and $c^0$ as follows.

\[
a^0_j = \begin{cases}
x_j-1 & \text{if } j< \gl_{1} \text{ or } j=n,\\
x_j & \text{otherwise},
\end{cases}
~~~~
~~~~
b^0_j = \begin{cases}
0 & \text{if } j< \gl_{1} \text{ or } j=n,\\
x_j & \text{otherwise}.
\end{cases}
\]

\[
c^0_j = \begin{cases}
x_j-1 & \text{if } j< \gl_{1},\\
0 & \text{if } j=n,\\
x_j & \text{otherwise},
\end{cases}
\]

\[
a^1_j = \begin{cases}
0 & \text{if } j< \gl_{1},\\
x_j-1 & \text{if } j=\gl_1,\\
x_j & \text{otherwise},
\end{cases}
~~~~~~~
b^1_j = \begin{cases}
0 & \text{if } j\leq  \gl_{1},\\
x_j & \text{otherwise}.
\end{cases}
\]

For $i=2,\dots ,k$ we set
\[
a^i_j = \begin{cases}
0 & \text{if } j\leq \gl_{i-1},\\
x_j-1 & \text{if } \gl_{i-1}<j\leq \gl_{i},\\
x_j & \text{otherwise},
\end{cases}
~~~~~~~
b^i_j = \begin{cases}
0 & \text{if } j\leq \gl_{i},\\
x_j & \text{otherwise}.
\end{cases}
\]

Note first that we have $a^0\geq c^0\geq b^0$ and $a^i\geq b^i$ for all indices $i=1,\dots ,k$. There exists a hyperedge $H^0$ such that $x\to a^0$, $x\to c^0$, and $x\to b^0$ are all $H^0$-moves. For $i=1,\dots ,k$ there exists a hyperedge $H^i\in\cH$ such that both $x\to a^i$ and $x\to b^i$ are $H^i$-moves. Furthermore, we have $m(a^0)=m(x)-1$, $m(c^0)=m(b^0)=0$, and $m(a^i)=m(b^i)=0$ for all $i=1,\dots ,k$. It follows that all positions in the set
\[
X=[c^0,a^0]\cup [b^0,c^0]\cup \bigcup_{i=1}^k [b^i,a^i]
\]
are long, since $x$ was assumed to be long, and they are all reachable from $x$ by a single move.

Let us note next that due to the definition of $\T_\cH$, condition (i) of Theorem \ref{t-2}, and Lemma \ref{l-l1} we have $\T_\cH(a^{i+1}) \geq \T_\cH(b^i)-1$ for $i=1,\dots ,k-1$. Furthermore $\T_\cH(a^0)= \T_\cH(x)-1$ by Lemma \ref{l-l1} and $\T_\cH(b^k)=m(x)$ because of condition (ii) of Theorem \ref{t-2}.

Therefore we can apply Lemma \ref{l-t-contiguity} to the pairs $(a^0,c^0)$, $(c^0,b^0)$, and  $(a^i,b^i)$, $i=1,\dots ,k$ and obtain
\[
\{\T_\cH(x')\mid x'\in X\} ~=~ [m(x),\T_\cH(x)-1].
\]
\qed

\end{document}